\documentstyle[leqno]{article}

\newtheorem{definition}{Definition}[section]
\newtheorem{lemma}{Lemma}[section]
\newtheorem{theorem}{Theorem}[section]
\newtheorem{proposition}{Proposition}[section]

\newtheorem{oss}{Remark}[section]

{\end{description}}

\newcommand{\binom}[2]%
{\left( \begin{array}{c} #1 \\ #2 \end{array} \right) }

{\hfill $\bullet$ \\}

\def\intem{\frac 1 {m\{Q^d\}}\int_{Q^d}}
\newcommand{\ima}{\hat{\bf i}}

\def\rif#1{{\rm (\ref{#1})}}

\def\reali{{\mbox{\bf R}}}
\def\com{{\mbox{\bf C}}}
\def\naturali{{\mbox{\bf N}}}

\newcommand{\forget}[1]{}

\begin{document}
\title{The spectral approximation
of multiplication operators via asymptotic (structured) linear
algebra}
\author{Stefano Serra-Capizzano
       \thanks{Dipartimento di Fisica e
       Matematica, Universit\`a dell'Insubria,
       Via Valleggio 11, 22100 Como (ITALY).
       Email:\ stefano.serrac@uninsubria.it; serra@mail.dm.unipi.it}
       }
\maketitle
\date{}

\begin{abstract}
A multiplication operator on a Hilbert space may  be approximated with finite sections by choosing an orthonormal
basis of the Hilbert space. Nonzero multiplication operators on $L^2$ spaces of functions are never compact and
then such approximations cannot converge in the norm topology. Instead, we consider how well the spectra
of the finite sections approximate the spectrum of the multiplication operator whose expression is simply
given by the essential range of the symbol (i.e. the multiplier). We discuss the case of real
orthogonal polynomial bases and the relations with the classical Fourier basis whose choice leads to well studied
Toeplitz case. The use of circulant approximations leads to constructive algorithms working for the
separable multivariate and matrix-valued cases as well.
\end{abstract}

\section{Introduction}\label{sec:1}

This note is in some sense a consequence of the intriguing MathSciNet Revue by Albrecht B\"ottcher
of a paper by Morrison \cite{morrison} and, of course, of the intriguing paper itself.
Briefly, if $\phi$ is a bounded function
defined on a compact set $K$ of ${\bf R}^d$, $d\ge 1$, consider the multiplication
operator $M[\phi]:L_w^2(K)\rightarrow L_w^2(K)$ defined as $M[\phi](f)=\phi f$, $w$ suitable weight
function. It is known that the spectrum
is given by the essential range of $\phi$: now suppose that we have only a finite number of coefficients
$\left(M_N[\phi]\right)_{i,j}=\langle M[\phi]e_j,e_i\rangle$, $i,j=0,\ldots,N-1$,
with $\{e_j\}$ denoting an orthonormal basis of $L_w^2$;
the question is about the reconstruction of the multiplier $\phi$ from the spectra of $M_N[\phi]$.
For reconstruction we mean the convergence of the finite sections spectra to the the essential range
of the symbol $\phi$. More in general, we are interested in understanding as much as possible about $\phi$
only using the entries of the matrices $M_N[\phi]$ for large but finite $N$.

Indeed the problem posed is a classical one (a beautiful historical account can be found in \cite{morrison}).
Here the idea is to discuss how the case of the choice of
a general real orthogonal basis on $K=[-1,1]$ can be reduced to the Fourier case and therefore to the Toeplitz
case (see \cite{BS,CN,KS} and references therein
for an encyclopedic coverage from three different angles)
and how the latter can be reduced to the circulant case. Circulants (see \cite{Da})
are normal matrices since they can all be diagonalized by the same unitary transform.
Further the transform is the celebrated discrete Fourier transform (DFT) for
which a stable and extremely efficient algorithm exists (the Fast Fourier Transform i.e. FFT,
see \cite{fft}). Therefore
the general case can be translated into a problem of (asymptotic) structured numerical linear algebra for
which an accurate solution can be determined with a low computational cost (here for low
cost we mean $O(N\log(N))$ arithmetic operations i.e. the asymptotic cost of a generic FFT). Moreover,
the restriction on the boundedness of $\phi$ can be suppressed and, more precisely, a related symbol
$\tilde \phi$ (more specifically $\tilde \phi(x)=\phi(x)w(x)\sqrt{1-x^2}$) has to be supposed just Lebesgue integrable: in this case,
the operator $M[\phi]$ can be unbounded and has to be defined on a different domain. Multidimensional block
generalizations (for multiplication operators having a matrix-valued multivariate function as multiplier) are also
available thanks to the rich theory built in the finite dimensional case in recent years.
The paper contains three more sections: Section \ref{sec:2} is devoted to linear algebra tools; in Section
\ref{sec:3} we discuss the solution to our problem and we give a brief account on separable multivariate
and matrix-valued generalizations; Section \ref{sec:4} is concerned with open questions and final remarks.

\section{Notation from asymptotic linear algebra}\label{sec:2}

First we introduce some notations and definitions concerning
general sequences of matrices. For any function $F$ defined on $\com$ and
for any matrix  $A_n$ of size $d_n$, with eigenvalues $\lambda_j(A_n)$ and singular values
$\sigma_j(A_n)$, $j=1,\ldots,d_n$, by the symbols $\Sigma_{\sigma}(F,A_n)$ and $\Sigma_{\lambda}(F,A_n)$ we
denote the means
\[
{\frac 1 {d_n} \sum_{j=1}^{d_n}
F[\sigma_j(A_n)]}, \ \ \
{\frac 1 {d_n} \sum_{j=1}^{d_n}
F[\lambda_j(A_n)]},
\]
and by the symbol $\|\cdot\|$ the spectral norm i.e. $\|X\|$ is the maximal singular value
of the matrix $X$ (see \cite{Bha}).
Furthermore $\|\cdot\|_p$ indicates the Schatten $p$ norms, $p\in [1,\infty)$  defined as
\[
\|A_n\|_p^p=\Sigma_{\sigma}(|\cdot|^p,A_n)\cdot d_n.
\]
The Schatten $\infty$ ($p=\infty$) norm is exactly the spectral norm (for a unified treatment of these
norms refer to the beautiful book by Bhatia \cite{Bha}).
Moreover, given a sequence $\{A_n\}$ of matrices of size $d_n$ with
$d_n<d_{n+1}$ and given a $\mu$-measurable function $g$ defined over a set $K$ equipped with
a $\sigma$ finite measure $\mu$, we say that $\{A_n\}$ is distributed as
$(g,K,\mu)$ in the sense of the singular values (in the sense of
the eigenvalues) if for any continuous $F$ with bounded support
the following limit relation holds
\begin{equation}\label{distribution:sv-eig}
\lim_{n\rightarrow \infty}\Sigma_{\sigma}(F,A_n)={1\over \mu(K)}\int_K F(|g|)d\mu, \ \ \
\left(\lim_{n\rightarrow \infty}\Sigma_{\lambda}(F,A_n)=
{1\over \mu(K)}\int_K F(g)d\mu\right).
\end{equation}
In this case we write in short $\{A_n\}\sim_{\sigma} (g,K,\mu)$
($\{A_n\}\sim_{\lambda} (g,K,\mu)$). An interesting connection between the notion of distribution
and the Schatten $p$ norms is given in the following Lemma.

\begin{lemma} \label{distr-schatten}
Assume that $\{A_n\}\sim_{\sigma} (g,K,\mu)$ and that
$\|B_n\|_p=o(d_n^{1/p})$, $A_n,B_n$ both of size $d_n$, and $p\in
[1,\infty]$. Then it holds
\begin{equation}
\label{tesi-distr-schatten}
 \{B_n\}\sim_{\sigma} (0,K,\mu)\ \ \ \ {\rm and} \ \ \ \
 \{A_n+B_n\}\sim_{\sigma} (g,K,\mu).
\end{equation}
Moreover, if all the involved sequences are Hermitian and $\{A_n\}\sim_{\lambda} (g,K,\mu)$, then
(\ref{tesi-distr-schatten}) holds true with $\sim_{\sigma}$ replaced by $\sim_{\lambda}$.
\end{lemma}
{\bf Proof.} The tools for the proof in the case of $p=2$ can be found in \cite{tyrtyshnikov}.
Here we treat the general case by using analogous ideas. For $p=\infty$ and $\|B_n\|=o(1)$ the proof is trivial
by standard perturbation arguments (see e.g. \cite{Bha,Wil}). Therefore we focus our attention on the case where
$p\in [1,\infty)$. Indeed, from the assumptions on $\{B_n\}$ with $p\in [1,\infty)$, for every $\epsilon>0$,
we have
\begin{eqnarray*}
C(n) & = & \|B_n\|_p^p=\sum_{j=1}^{d_n} \sigma_j^p(B_n)  \\
  & \ge &  \sum_{\sigma_j(B_n)>\epsilon}\sigma_j^p(B_n) \\
  & \ge & \sum_{\sigma_j(B_n)>\epsilon} \epsilon^p \\
  & = &   \epsilon^p \#\{\sigma_j(B_n)>\epsilon\}
\end{eqnarray*}
with $C(n)=o(d_n)$. Therefore the cardinality of the singular
values bigger than $\epsilon$ is bounded from above by
$C(n)/\epsilon^p=o(d_n)$. Since $\epsilon>0$ is arbitrary, by
direct check, it follows that $\{B_n\}\sim_{\sigma} (0,K,\mu)$.
Furthermore, by exploiting the singular values decomposition of
$B_n$, we can write $B_n$ as $L_n(\epsilon)$ and $R_n(\epsilon)$
where $\|L_n(\epsilon)\|_\infty \le \epsilon$ and the
rank$(R_n(\epsilon)\le C(n)/\epsilon^p=o(d_n)$. More precisely, in
the previous lines we have proved that the cardinality of the
singular values of $B_n$ bigger than $\epsilon$ is bounded from
above by $C(n)/\epsilon^p=o(d_n)$. Now from the SVD decomposition
(see e.g. \cite{Bha}) there exist $U_n$ and $V_n$ unitary matrices
and $D_n$ diagonal matrix (containing the singular values of $B_n$
sorted non decreasingly) such that
\[
B_n = U_n D_n V_n.
\]
At this moment take $D_n(>)$ the matrix containing all the entries
bigger than $\epsilon$ of $D_n$ (in the same position as $D_n$)
and $D_n(<)$ the matrix containing all the entries at most equal
to $\epsilon$ of $D_n$ (in the same position as $D_n$). Therefore
$D_n=D_n(>)+D_n(<)$ with
\[
\|D_n(<)\|\le \epsilon,\ \ \ \ \ {\rm rank}(D_n(>))\le
C(n)/\epsilon^p=o(d_n).
\]
Finally since $U_n$ and $V_n$ are unitary we have
\[
\|U_nD_n(<)V_n\|=\|D_n(<)\|\le \epsilon,\ \ \ \ \ {\rm
rank}(U_nD_n(>)V_n)={\rm rank}(D_n(>))\le C(n)/\epsilon^p=o(d_n),
\]
and $B_n=U_nD_n(<)V_n+U_nD_n(>)V_n$. The statement is proven by
putting $L_n(\epsilon)=U_nD_n(<)V_n$ and
$R_n(\epsilon)=U_nD_n(>)V_n$

Consequently, by using e.g. Proposition 2.3 and Remark 2.1 in
\cite{algebra}, from the hypothesis $\{A_n\}\sim_{\sigma}
(g,K,\mu)$ we deduce $\{A_n+B_n\}\sim_{\sigma} (g,K,\mu)$. The
case of the eigenvalues for Hermitian matrices $A_n$ and $B_n$ is
identical and it is not repeated here. \ \hfill $\bullet$
\ \\

\subsection{How to use spectral distribution}\label{sec:2-1}

We show how the notion of distribution can be used for the reconstruction of the symbol when the
eigenvalues (or singular values) are known. More precisely, the subsequent Theorem \ref{th-relations-def}
demonstrates that $\{A_n\}\sim_{\lambda} (g,K,\mu)$ (or $\{A_n\}\sim_{\sigma} (g,K,\mu)$) and the knowledge
of the eigenvalues of $\{A_n\}$ (or singular values of $\{A_n\}$) imply that many facts on the symbol
$g$ can be constructively recovered.

\begin{definition}\label{def-range}
Given the $\mu$ measurable function $g$ defined on $K$ with $\mu$ being a $\sigma$ finite measure
supported on $K$, the (essential) range of $g$ is given by the points $s\in \bf C$
such that, for every $\epsilon>0$, the measure
of the set $\{t\in D:\ g(t)\in D(s,\epsilon)\}$ is positive with
$D(s,\epsilon)=\{z\in \bf{C}: |z-s|<\epsilon\}$. The function $g$ is (essentially) bounded if
its essential range is bounded. Finally, if $g$ is real-valued then the (essential) supremum
is defined as the supremum of its range  and the (essential) infimum
is defined as the infimum of its range.
\end{definition}

\begin{definition}\label{def-cluster}
A sequence $\{A_n\}$ ($A_n$ of size $d_n$) is properly (or strongly)
clustered at $s \in \bf C$ in the eigenvalue sense,
if for any $\epsilon>0$ the number of the eigenvalues of $A_n$ not
belonging to $D(s,\epsilon)=\{z\in \bf{C}: |z-s|<\epsilon\}$
can be bounded by a pure constant $q_\epsilon$ possibly depending on $\epsilon$ but not on $n$.
Of course if every $A_n$ has, at least definitely, only real eigenvalues, then
$s$ has to be real and the disk $D(s,\epsilon)$ reduces to the interval
$(s-\epsilon,s+\epsilon)$.
Furthermore, a sequence $\{A_n\}$ ($A_n$ of size $d_n$) is properly (or strongly)
clustered at the nonempty closed set $S \subset \bf C$ in the eigenvalue sense
if for any $\epsilon>0$ the number of the eigenvalues of $A_n$ not
belonging to $D(S,\epsilon)=\displaystyle \bigcup_{s\in S}D(s,\epsilon)$
can be bounded by a pure constant $q_\epsilon$ possibly depending on $\epsilon$ but not on $n$ and
if every $A_n$ has, at least definitely, only real eigenvalues, then
$S$ has to be a nonempty closed subset of  $\bf R$.
The term ``properly (or strongly)'' is replaced by ``weakly'' if $q_\epsilon$ is a possibly unbounded
function of $n$ with $q_\epsilon(n)=o(d_n)$
(i.e. $\displaystyle \lim_{n\rightarrow\infty} {q_\epsilon(n)\over d_n}=0$).
Finally, the above notions are in the singular value sense if the term ``eigenvalue'' is replaced by
``singular value'': of course $s$ has to be a real nonnegative number and $S$ has to be
a subset of nonnegative numbers.
\end{definition}

\begin{definition}\label{def-attractors}
A sequence $\{A_n\}$ ($A_n$ of size $d_n$ and with spectrum $\Sigma_n$)
is strongly attracted by $s \in \bf C$ if
\[
\lim_{n\rightarrow \infty}{\rm dist}(s,\Sigma_n)=0
\]
where ${\rm dist}(X,Y)$ is the usual Euclidean distance between two subsets $X$ and $Y$ of the complex plane.
Furthermore, let us order the eigenvalues according to its distance from $s$ i.e.
$|\lambda_1(A_n)-s|\le |\lambda_2(A_n)-s| \le \cdots \le |\lambda_{d_n}(A_n)-s|$. We say that the attraction is of
order $r(s)\in \bf N$, $r(s)\ge 1$, fixed number independent of $n$, if
\[
\lim_{n\rightarrow \infty}|\lambda_{r(s)}(A_n)-s|=0,\ \ \ \
\liminf_{n\rightarrow \infty}|\lambda_{r(s)+1}(A_n)-s|>0.
\]
The attraction is of order $r(s)=\infty$ if
\[
\lim_{n\rightarrow \infty}|\lambda_{j}(A_n)-s|=0
\]
for every fixed $j$ independent of $n$. Finally, the term ``strong or strongly'' is replaced
by ``weak or weakly'' if every symbol $\lim$ is replaced by $\liminf$.
Finally, the above notions are in the singular value sense if the term ``eigenvalue'' is replaced by
``singular value'', $\Sigma_n$ is replaced by the set of the singular values,
and, of course, the value $s$ is a real nonnegative number.
\end{definition}

Notice that writing $\{A_n\}\sim_{\lambda} (g,K,\mu)$ with $g$ constant function equal
to $s\in \bf C$ is equivalent to write that  $\{A_n\}$ is weakly clustered at $s$
in the eigenvalue sense. Analogously, writing $\{A_n\}\sim_{\sigma} (g,K,\mu)$ with $g$ constant function equal
to $s\in \bf R$, $s\ge 0$, is equivalent to write that  $\{A_n\}$ is weakly clustered at $s$
in the singular value sense.

The notions previously introduced are intimately related as emphasized in the subsequent theorem which is
explicitly given only for the eigenvalues (the singular value version is obvious and is shortly sketched).

\begin{theorem}\label{th-relations-def}
Let $\{A_n\}$ be a matrix sequence with $A_n$ having size $d_n$
and let $g$ be a $\mu$-measurable function defined on $K$ with
$\mu$ being $\sigma$ finite measure supported on $K$. Consider the
following statements:
\begin{description}
\item[a)] $\{A_n\}\sim_{\lambda} (g,K,\mu)$;
\item[b)] the (essential) range of $g$ is a weak cluster for $\{A_n\}$ in the eigenvalue sense;
\item[c)] the (essential) range of $g$ strongly attracts the eigenvalues of $\{A_n\}$;
\item[d)] any point $s$ of the (essential) range of $g$ strongly attracts the eigenvalues of $\{A_n\}$
with order $r(s)=\infty$.
\item[e)] given $s\in \bf C$, $\epsilon>0$, if the cardinality of the eigenvalues of $A_n$ belonging to
$D(s,\epsilon)$ divided by $d_n$ tends to a positive value, then $s$ belongs to the
(essential) range of $g$ within an error of at most $\epsilon$;
\item[f)] given $s\in \bf C$, $\epsilon>0$, if the cardinality of the eigenvalues of $A_n$ belonging to
$D(s,\epsilon)$ divided by $d_n$ tends to a zero, then $s$ cannot belong to the
(essential) range of $g$.
\end{description}
Then {\bf a)} implies {\bf b)}, {\bf c)}, {\bf d)}, {\bf e)}, and {\bf f)}.
Finally, the above implications hold in the singular value sense if the term ``eigenvalue'' is replaced by
``singular value'', $g$ is replaced by $|g|$,
and of course the value $s$ is a real nonnegative number.
\end{theorem}
{\bf Proof.}
The first three implications are proven in Theorem 2.7 of \cite{Salmeria}. For the other two
see e.g. Section 4 in {\rm \cite{Skoro}}.
\ \hfill $\bullet$
\ \\

In the rest of the paper, with regard to relationships (\ref{distribution:sv-eig}), the symbol $\mu$
is suppressed for the cases under study (Toeplitz sequences, Generalized
Locally Toeplitz sequences, Circulants etc.) since the measure will
always  coincide with the standard Lebesgue measure on $\reali^d$ for
some positive integer $d$.

\subsection{Toeplitz matrix sequences}\label{subs:2-1}

Let $m\{\cdot\}$ be the Lebesgue measure on ${\bf R}^d$ for some $d$ and
let $f$ be a $d$ variate complex-valued (Lebesgue) integrable function,
defined over the hypercube $Q^d$, with $Q=(-\pi,\pi)$ and $d\ge 1$.
From the Fourier coefficients of $f$
\begin{equation}
\label{defcoeff}
f_{j}=\intem f(s) {\rm exp}({-\ima (j,s)})\, ds,\qquad \ima^2=-1,\quad
j=(j_1,\ldots,j_d)\in {\bf Z}^d
\end{equation}
with $(j,s)=\sum_{k=1}^d j_k s_k$, $n=(n_1,\ldots,n_d)$ and
$N(n)=n_1\cdots n_d$, we can build the sequence of Toeplitz matrices
$\{T_n(f)\}$, where $T_n(f)=$$\{f_{j-i}\}_{i,j={\bf 1}^T}^n$
$\in {\bf C}^{N(n)\times N(n)}$,
${\bf 1}^T=(1,\ldots,1)\in \naturali^d$ is said to be the Toeplitz matrix of order
$n$ generated by $f$. Furthermore, throughout the
paper when we write $n\rightarrow \infty$ with $n=(n_1,\ldots,n_d)$ being
a multi-index, we mean that $\min_{1\le j\le d} n_j\rightarrow \infty$.\par
The asymptotic distribution of eigen and singular values of a sequence of
Toeplitz matrices  has been thoroughly studied in the last century (for
example see \cite{BS}\ and the references reported therein).
Here we report a famous Theorem  of Szeg\"o \cite{szego}, which we state in the
Tyrtyshnikov and Zamarashkin version \cite{tyrtL1}:
\begin{theorem}
\label{teoszego-tyr}
If $f$ is integrable over $Q^d$, and if $\{T_n(f)\}$  is the sequence of
Toeplitz matrices generated by $f$, then it holds
\begin{equation}
\label{tesityr}
 \{T_n(f)\}\sim_{\sigma} (f,Q^d).
\end{equation}
Moreover, if $f$ is also real-valued, then each matrix $T_n(f)$
is Hermitian and
\begin{equation}
\label{tesiszego}
\{T_n(f)\}\sim_{\lambda} (f,Q^d).
\end{equation}
\end{theorem}

This result has been generalized to the case where $f$ is
matrix-valued (see, for example, \cite{Tillinota,Smarko1} and Subsection \ref{subs:3-3}) so that the
matrices $T_n(f)$ have multilevel block Toeplitz structure and to the case where the test
functions $F$ have not bounded support (see \cite{Spal} and references therein).
\par
If $f$ is not real-valued, then $T_n(f)$ is not Hermitian in general: consequently, the distribution of
eigenvalues is more involved and \rif{tesiszego}\ cannot be extended in the natural way (see \cite{tilliLA}).
A very elegant geometric based result is due to Tilli \cite{tillicomplex} and the conclusion is surprisingly
simple: a Toeplitz sequence with bounded symbol $f$ will have a canonical eigenvalue distribution in the sense
of (\ref{distribution:sv-eig}), if the complement of the range of $f$ is connected and the range has empty interior.
This makes clear that regularity plays no role and this explain why
this result was not found for many years: researchers were in the wrong direction looking at regularity
assumptions on the symbol. The same misunderstanding occurred, in minor proportions, for the conditioning of a
Toeplitz matrix generated by a weakly sectorial symbol \cite{BG}: again it is a geometric phenomenon that
describes the asymptotic behavior of the conditioning and not a regularity property of the symbol.
Take $f(s)=(2-2\cos(s))^{10}$. Then the minimal eigenvalues of the single-level $T_n(f)$ tends to $0$
(the infimum of $f$) monotonically and with asymptotic speed dictated by $n^{-20}$ (notice that $20$ is the order
of the unique zero of $f$). Exactly the same behavior
is proven (with a different constant \cite{Se1,BG}) if $f(s)=(2-2\cos(s))^{10}h(s)$ where $h(s)$ is any
real-valued $L^\infty$ function with positive infimum: indeed the result is a consequence of how the
essential range of the nonnegative symbol $f$ ``touches'' $0$ from above and the fact that $f(s)$ is
infinitely differentiable, as in the case of $h(s)=1$, or is discontinuous almost everywhere (a.e.) does
play any role.

\subsection{GLT matrix sequences}\label{subs:2-2}

For the subsequent analysis, it is convenient to introduce the class of Generalized Locally Toeplitz (GLT)
sequences that represents at the same time a generalization of Toeplitz sequences and of matrix sequences
approximating variable coefficient (differential) operators \cite{multiloc}.
More in detail, the class of GLT sequences can be essentially viewed as a topological closure,
both in the matrix side and in the ``symbol'' side, of linear combinations of products of Toeplitz sequences
and diagonal sampling matrix sequences: a sampling matrix (of level $1$) $D_n(a)$ of size $n$ is the diagonal
matrix containing as $j$-th diagonal element $a(jh)$, $h$ a mesh parameter, $a$ smooth enough. For our
purposes we do not need to introduce the (quite long and involved) definition of this class for which we refer
to \cite{multiloc,sccm}. We just recall the main properties especially those which are of interest
for our problem.

\begin{description}
\item[A.]
Any GLT sequence $\{A_n\}$ is uniquely associated to a measurable symbol
$\kappa(x,s)$, $x\in \Omega$ Peano-Jordan measurable
set of ${\bf R}^d$ (space domain), $s\in Q^d$ (Fourier domain), $D=\Omega\times Q^d$:
we write $\{A_n\}\sim_{\rm GLT} \kappa$ and we have
$\{A_n\}\sim_{\sigma} (\kappa,D)$ and $\{A_n\}\sim_{\lambda} (\kappa,D)$
if $A_n$ Hermitian at least for $n$ large enough.
\item[B.]
Every Toeplitz sequence generated by $f(s)$ in the sense of (\ref{defcoeff}) is a GLT sequence with
$\kappa(x,s)=f(s)$ (Szeg\"o-Tyrtyshnikov theory).
\item[C.]
Every sequence which is distributed as the zero function in the sense of (\ref{distribution:sv-eig})
for the singular value is a GLT sequence with $\kappa(x,s)=0$.
\item[D.]
Every Finite Difference (FD) and Finite Element (FE) equi-spaced approximations of
constant coefficient PDEs on square regions
(any boundary condition) is a GLT sequence with $\kappa(x,s)=p(s)$ for some trigonometric
polynomial $p$  (Fourier Analysis);
\item[E.]
Any FD, FD discretization of general variable coefficient (system of) PDEs over $\Omega$
is a GLT sequence. In that case $\kappa(x,s)$ is easily identified (Generalized Fourier Analysis):
$\kappa(x,s)$ is the {\em principal symbol} with obvious changes
of the Kohn-Nirenberg and H\"ormander theory for {\em Pseudo-Differential} operators.
\end{description}

The GLT sequences form a $*$--algebra. More precisely, the GLT sequences are stable under
linear combinations, product, pseudo-inversion, and adjoint. In fact, if
$\{A_n\}\sim_{\rm GLT}\kappa_A$ and $\{B_n\}\sim_{\rm GLT}\kappa_B$, then we observe stability under
\begin{description}
\item[F.] linear combinations i.e. $\{\alpha A_n+\beta B_n\}\sim_{\rm GLT}
\alpha\kappa_A+\beta \kappa_B$;
\item[G.] product i.e. $\{A_nB_n\}\sim_{\rm GLT} \kappa_A\kappa_B$;
\item[H.] (pseudo)-inversion i.e. $\{A_n^+\}\sim_{\rm GLT} \kappa_A^{-1}$
provided that $\{A_n\}$ is invertible (invertible elements are those such that
the symbol vanishes at most on a set of zero Lebesgue measure=sparsely vanishing).
\item[I.] adjoint i.e. $\{A_n\}\sim_{\rm GLT} \kappa_A$  is equivalent to $\{A_n^*\}\sim_{\rm GLT}
\kappa_A^*$.
\end{description}

In the following, we will use essentially properties {\bf A.}, {\bf B.}, {\bf C.}, and the structure
of algebra of the GLT class.

\section{Identification of the multiplier}\label{sec:3}

The section is divided in three parts. In the first we discuss in detail the solution to our problem
in one dimension: as already mentioned, it turns out that the boundedness of $\phi$ is not necessary
and only the Lebesgue integrability of a related symbol $\tilde \phi$ is crucial. Subsections \ref{subs:3-2}
and \ref{subs:3-3} are devoted to sketch the solution for multivariate and matrix-valued multipliers
in the case of separable weight functions. Instead of giving all the details, we will emphasize what is new
in the derivation and the surprise is that the multivariate matrix-valued problem does not pose essentially
more difficulties than the scalar case in one dimension (except, may be, for the notations).

\subsection{The problem in $1$ dimension}\label{sec:3-1}

Let $\phi$ be a bounded function defined on $[-1,1]$ and let us consider the multiplication
operator $M[\phi]:L_w^2([-1,1])\rightarrow L_w^2([-1,1])$ defined as $M[\phi](f)=\phi f$, $w$ suitable weight
function with the usual assumptions. It is known that the spectrum
is given by the essential range of $\phi$ (see e.g. \cite{morrison}):
now suppose that we have only a finite number of coefficients
\begin{equation}\label{form-m-phi}
M_n[\phi]=\left(\langle M[\phi]e_j,e_i\rangle\right)_{i,j=0}^{n-1},
\end{equation}
with $\{e_j\}$ denoting an orthonormal basis of $L_w^2$. The question concerns the reconstruction of the
multiplier $\phi$ from the spectrum of $M_n[\phi]$ in the sense discussed in Subsection \ref{sec:2-1}.

Consider first the case of the Chebyshev weight $w(x)=(1-x^2)^{-1/2}$ of first kind and of its basis
$e_j(x)=\cos(j\arccos(x))$. Then
\[
\left(M_n[\phi]\right)_{i,j}=\langle M[\phi]e_j,e_i\rangle=\int_{-1}^1\phi(x)e_j(x)\overline{e_i}(x)\, dx=
\int_{Q}\phi(\cos(s))\cos(js)\cos(is)\, ds,\ \ Q=(-\pi,\pi).
\]
As a consequence, the matrix $M_n[\phi]$ can be expressed in terms of the Fourier coefficients $f_k$ of the
function $f(s)=\pi \phi(\cos(s))$ in the sense of (\ref{defcoeff}) and then
$\left(M_n[\phi]\right)_{i,j}={1\over 2}(f_{i-j}+f_{j-i}+f_{i+j}+f_{-i-j}).$
Taking into account that $f(s)$ is even we directly see that $f_{k}=f_{-k}$ for ever $k\in \bf Z$ and therefore
\[
\left(M_n[\phi]\right)_{i,j}=f_{|i-j|}+f_{|i+j|},
\]
i.e.
\begin{equation}\label{toep-hankel}
M_n[\phi]=T_n(f)+H_n(f).
\end{equation}
Here the matrix $H_n(f)=\left(f_{|i+j|}\right)_{i,j=0}^{n-1}$ is of Hankel type since its entries are
constant along the anti-diagonals. Moreover, from \cite{FasTilli} we know that the Hankel sequence
$\{H_n(f)\}$ is distributed as the zero function over $Q$ in the sense
of (\ref{distribution:sv-eig}): $\{H_n(f)\}$ is indeed a GLT sequence \cite{multiloc} with symbol equal to zero
(see item {\bf C.}) and therefore, since it is bounded in spectral norm,
both relationships in (\ref{distribution:sv-eig}) hold with $g=0$ (the singular value part is contained
in item {\bf A.} and {\bf C.} and the eigenvalue part follows as in Theorem 1.2 of \cite{cluster}).
Consequently, the singular value distribution of $\{M_n[\phi]\}$ is decided by the Toeplitz part $\{T_n(f)\}$
and, if $\phi$ is real-valued, the same is true for the eigenvalue distribution too: this can
be seen directly by using Tyrtyshnikov perturbation arguments \cite{tyrtyshnikov} or, from a more abstract
viewpoint, because the GLT class is an algebra that is by {\bf A.} and {\bf F.}
(see also \cite{Tilliloc}). Therefore the symbol
of $\{M_n[\phi]\}=\{T_n(f)\}+\{H_n(f)\}$ is equal to the one of $\{T_n(f)\}$ i.e. $f$ plus that of
$\{H_n(f)\}$ which is zero.

As a consequence, if the multiplier $\phi$ is real-valued, then we can reconstruct, approximately, $\phi$ from the
eigenvalues of $M_n[\phi]$. In the general case, the desired result depends on the geometric structure of the
range of $\phi$ and on the Hankel correction: the eigenvalues can be dramatically sensitive even to $1$ rank
corrections (see e.g. \cite{Wil}). This pathological behavior of the eigenvalues has also good side effects
because the effective procedure that can be designed (see Subsection \ref{subs:3-1-2}) depends exactly on the
existence of close sequences whose spectral behavior is substantially more regular than Toeplitz
sequences (see Remark \ref{circulants-vs-rest-of-the-world}).

The case of the Chebyshev weight of second kind is also very simple to handle thanks to the explicit expression
of its orthogonal basis elements after the usual change of variable $x=\cos(s)$. Indeed we have
$w(x)=(1-x^2)^{1/2}$ and $e_j(x)=\sin((j+1)\arccos(x))/\sin(\arccos(x))$. Then
\[
\left(M_n[\phi]\right)_{i,j}=\langle M[\phi]e_j,e_i\rangle=\int_{-1}^1\phi(x)e_j(x)\overline{e_i}(x)\, dx=
\int_{Q}\phi(\cos(s))\sin((j+1)s)\sin((i+1)s)\, ds,\ \ Q=(-\pi,\pi).
\]
From the latter we infer $\left(M_n[\phi]\right)_{i,j}=f_{|i-j|}+f_{|i+j+2|}$ with $f_k$ Fourier coefficients
of $f(s)=\pi \phi(\cos(s))$ and then $M_n[\phi]=T_n(f)+\tilde H_n(f)$ with $\tilde H_n(f)$ being the principal
sub-matrix (of size $n$) made by the last $n$ rows and columns of $H_{n+1}(f)$ and
$H_{n}(f)$ as in (\ref{toep-hankel}).
Therefore a simple interlace argument (see e.g. \cite{Bha})
shows that the corresponding Hankel sequence is distributed as the zero function over $Q$ and then
(see \cite{multiloc}) $\{M_n[\phi]\}=\{T_n(f)\}+\{H_n(f)\}$ has the same GLT symbol as $\{T_n(f)\}$
i.e. $f$ and the conclusion is as before.

In fact, the above analysis can be generalized using purely linear algebra tools but the result itself is known
already thanks to Szeg\"o (see \cite{szego-pol}).
For every choice of the weight function the symbol of $\{M_n[\phi]\}$ is always
$f(s)=\pi \phi(\cos(s))$ which is independent of the weight function $w$. In other words, the finite sections
of $M[\phi]$ with orthogonal polynomials always give more attention to the endpoints of the original
interval $-1$ and $1$ and less attention on the central part of the domain. That behavior is also important
for the success of many associated numerical methods such as Gaussian quadrature formulae and
interpolations schemes at the zeros of orthogonal polynomials.

However, let us give a short look to a sketch of a linear algebra derivation.

\begin{proposition}\label{prop:sec:3-1}
Consider a general weight $w$ with the usual restrictions (nonnegative, with support coinciding with
$[-1,1]$, with finite Lebesgue integral). Let $e_j$ be the $j$-th orthogonal polynomial. Then the following
facts hold:
\begin{enumerate}
\item $e_j(x)=\sum_{i=0}^{j} a_i c_i(x)$, $a_j\neq 0$, $c_i$ $i$-th Chebyshev polynomial of first kind;
\item $E_{n-1}(x)=L_n F_{n-1}(x)$, $L_n$ lower triangular invertible matrix,
$E_{n-1}(x)$ $n$-dimensional vector whose $i$-th position, $i=0,\ldots,n-1$, is the given
by $e_i(x)$ and $F_{n-1}(x)$ $n$-dimensional vector whose $i$-th position, $i=0,\ldots,n-1$, is the given
by $c_i(x)$;
\item $M_n[\phi]=\int_{-1}^1\phi(x)w(x) E_{n-1}(x)E_{n-1}^H(x)\, dx$;
\item $M_n[\phi]=L_n\cdot \left[\int_{-1}^1\phi(x)w(x) F_{n-1}(x)F_{n-1}^H(x)\, dx\right] \cdot L_n^H$;
\item $M_n[\phi]=L_n\cdot \tilde M_n[\tilde \phi] \cdot L_n^H$, with $\tilde \phi(x)=\phi(x) w(x) \sqrt{1-x^2}$
and $\tilde M_n[\tilde \phi]$ being the $n$-th finite section of $M[\tilde \phi]$ in the case of the
Chebyshev weight of first kind;
\item $\{\tilde M_n[\tilde \phi]\}\sim_\sigma (\tilde f,Q)$,
$\tilde f(s)=\pi\tilde \phi(\cos(s))=\pi\phi(\cos(s))w(\cos(s))\sin(s)$, $Q=(-\pi,\pi)$;
\item $\{L_n\}$ is a GLT sequence with symbol $\sqrt{g(s)}$ and $\{L_n^H L_n\}$
is a GLT sequence with symbol $g(s)=1/[w(\cos(s))\sin(s)]$;
\item $\{L_n^H L_n\}\sim_{\lambda,\sigma} (g,Q)$, $g(s)=1/[w(\cos(s))\sin(s)]$, $Q=(-\pi,\pi)$;
\item $\{ M_n[ \phi]\}$ is a GLT sequence with weight $f(s)=\pi \phi(\cos(s))$
\item $\{ M_n[ \phi]\}\sim_\sigma (f,Q)$, $f(s)=\pi \phi(\cos(s))$.
\end{enumerate}
\end{proposition}
{\bf Discussion on the proof.}\
{\bf Item 1.} is obvious since every $e_j$ has degree $j$ and the second item is again obvious
since $e_j$ has exactly degree $j$. {\bf Item 3.} is a compact rewriting, directly in matrix form,
of (\ref{form-m-phi}). {\bf Item 4.} follows from {\bf Item 2.} and {\bf Item 3.} taking into account
the linearity of the integral and that $L_n$ does not depend on $x$.
We now recall that $F_{n-1}(x)$ contains the Chebyshev basis of first kind: therefore, in order to interpret
the scalar product as the one induced by the Chebyshev weight of first kind, the related multiplier has
to be seen as $\tilde \phi(x)=\phi(x) w(x) \sqrt{1-x^2}$ and {\bf Item 5.} is proved.
Further, after the usual change of variable $x=\cos(s)$, the matrix $\tilde M_n[\tilde \phi]$ can be
written as $T_n(\tilde f)$ plus $H_n(\tilde f)$. Moreover, $w\in L^1[-1,1]$ and therefore $\tilde f\in L^1(Q)$,
$Q=(-\pi,\pi)$. Finally by Theorem \ref{teoszego-tyr} (which holds for $L^1$ functions) and by
\cite{FasTilli}, we know that the Toeplitz part is distributed as $f$ and the Hankel part as zero,
respectively (also {\bf Item A.}, {\bf Item B.}, and {\bf Item C.}).
This is enough by {\bf Item F.} (see Theorem 4.5 and Section 5 in \cite{multiloc} for more details)
for deducing that $\{\tilde M_n[\tilde \phi]\}\sim_\sigma (\tilde f,Q)$
and {\bf Item 6.} is proven. The first part of {\bf Item 7.} is just a technical calculation (see \cite{nota})
and the second part of the same item and {\bf Item 8.} follow from the first part and from {\bf Item G.} and
{\bf Item I.} (see also Theorem 4.5 in \cite{multiloc}).
Finally, by {\bf Item 5.} and {\bf Item 7.}, $\{ M_n[ \phi]\}$
is a product of two GLT sequences, $\{\tilde M_n[\tilde \phi]\}$ and $\{\tilde L_n^H L_n\}$ with symbols
$\pi\phi(\cos(s))w(\cos(s))\sin(s)$ and $1/[w(\cos(s))\sin(s)]$, respectively. Therefore, due to the structure
of algebra of GLT sequences (again Theorem 5.8 in \cite{multiloc} i.e. {\bf Item G.}),
$\{ M_n[ \phi]\}$ is a GLT sequence with symbol
$\pi\phi(\cos(s))=\pi\phi(\cos(s))w(\cos(s))\sin(s)/[w(\cos(s))\sin(s)]$ ({\bf Item 9.}) and,
finally, $\{M_n[ \phi]\}\sim_\sigma (f,Q)$, $f(s)=\pi \phi(\cos(s))$ that is {\bf Item 10.},
again by Theorem 4.5 in \cite{multiloc} i.e. {\bf Item A.}\
\hfill $\bullet$
\ \\

It is clear that the above proof is really compressed and that some relevant details have to be found
in other papers (in particular for {\bf Item 7.}).
The point was to show that from purely linear algebra reasonings it is possible to treat
this kind of problems and sometimes obtaining in a simpler way more general information: see e.g.
\cite{ku-ser} where the analysis of the zero distribution of orthogonal polynomials with varying coefficients
is made by employing GLT arguments, without any regularity assumption except for the Lebesgue measurability.
We emphasize in addition that the algorithm
in the next subsection depends only on {\bf Items 4., 5.}, and {\bf 6.} for which we have given a detailed proof
and that item {\bf Item 6.} is indeed valid as long as $\tilde f \in L^1(Q)$. We observe that the latter means
that the assumption on the boundedness of the multiplier $\phi$ is not necessary and can be dropped.
More specifically, we can allow $\phi$ to be just Lebesgue integrable if we have
$w(\cos(s))\sin(s)\in L^\infty(Q)$: we
already encountered examples in this direction and namely the Chebyshev weight of first kind
for which $w(\cos(s))\sin(s)=1$ and that of second kind for which $w(\cos(s))\sin(s)=\sin^2(s)$.
It is clear that there exists a large class of weights of this type.

\subsubsection{Circulant approximation}\label{subs:3-1-1}

We start by describing the circulant class with special attention to its approximation
properties with respect to Toeplitz matrix sequences. The algebra of circulant matrices
is a subclass of Toeplitz matrices to which it is not possible to attribute a symbol
in the sense of (\ref{defcoeff}) with exception for the identity and for the null matrix.
In the one-level case (the one discussed so far in this section), they share
the algebraic property that every row is the forward circular one-step shift of the previous
row and where also the notion of ``previous'' has to be intended in a circular way: more precisely,
the first row can be seen as the forward circular one-step shift of the last row as it is
clear from equation (\ref{circ-ex}). The latter
nice algebraic feature translates in many properties related to circular convolutions. Here we only
point out another important characterization in a spectral sense. Every circulant matrix of size $n$
can diagonalized by the (unitary) discrete Fourier matrix. This means that $A_n$ is circulant
if and only if $A_n=F_n D F_n^*$ where $D$ is a complex diagonal matrix,
\[
F_n= \left({\displaystyle {1\over \sqrt{n}}e^{-2\pi{\bf  i}jk/n
}}\right), \ \ \ k,j=0,\ldots,n-1,
\]
is the Fourier matrix of size $n$ and $X^*$ denotes the complex transpose of $X$.
Moreover, the diagonal matrix $D$ has $j$-th entry given by $p_n(x_j^{(n)})$
with $x_j^{(n)}=2\pi j/n$, $j=0,\ldots,n-1$, $p_n(z)=\sum_{k=0}^{n-1}a_k z^k$, $a_0,\ldots,a_{n-1}$ being
the entry of the first column $c[1]$ of $A_n=$circ$(a)$ i.e.
\begin{equation}\label{circ-ex}
A_n=\left[ \begin{array}{ccccc}
    a_{0} & a_{n-1} & \cdots & a_2 & a_1 \\
    a_{1} & a_{0}  & a_{n-1} & \cdots & a_2 \\
    a_2    & \ddots &\ddots& \ddots & \vdots\\
    \vdots &\ddots &\ddots &\ddots& a_{n-1} \\
    a_{n-1}  &\cdots  &    a_2  & a_{1} & a_{0}
     \end{array} \right].
\end{equation}

Notice that the above eigenvalue formula has also
an important computational counterpart since the vector $d$ containing the diagonal entries $D$
is equal to $F_n^* c[1]$ and $F_n^*=PF_n$, with $P$ flip-type permutation matrix. As a consequence, the
spectral decomposition of any circulant matrix can be recovered in $O(n\log(n))$ complex
operations via the celebrated FFT (see \cite{fft}).
We now recall some connections between circulants and (one-level) Toeplitz matrix sequences
associated to a symbol.

\begin{definition}\label{approx-circ}
Let $C_n$ be the algebra of circulant matrices and let $T_n(f)$ be
a single level Toeplitz matrix associated to the symbol $f$. Then the following
definitions hold.
\begin{itemize}
\item The Strang preconditioner $N_n(f)$ associated to $T_n(f)$ is the circulant matrix
obtained from $T_n(f)$ by copying the first $[n/2]$ central diagonals with $[x]$ denoting
the rounding of $x$. In other words, the $j$-th entry of first column $c[1]$ of $N_n(f)$,
$j=0,\ldots, [n/2]-1$, is exactly the $j$-th Fourier coefficient $a_j$ of $f$.
\item The optimal preconditioner $C_n(f)={\rm Opt}(T_n(f))$ is the unique solution of the minimization
problem
\begin{equation}\label{minsuper}
\min_{X\in C_n} \|A-X\|_{\rm F},\ \ \ \ A=T_n(f),
\end{equation}
with $\|\cdot\|_F$ denoting the Frobenius norm i.e. the Euclidean norm
of the singular value vector (Schatten $p$ norm with $p=2$)
or, equivalently, the Euclidean norm of $n^2$-sized
vector obtained by putting in a unique vector all the columns of the argument.
\end{itemize}
\end{definition}

Some remarks are in order. The existence and uniqueness of the Strang or natural preconditioner
are implicit in the definition itself which clearly indicates an explicit cost-free
expression. The existence and uniqueness of the optimal preconditioner (see e.g. \cite{CN})
follow from the strict convexity of the Frobenius norm
that implies the existence and uniqueness of the minimizer from a given convex closed set. We are
in a finite dimensional setting and clearly the linear space of the circulants $C_n$ is closed and convex.

Finally, the optimal approximation admits an easy to derive and very interesting representation
since
\begin{equation}\label{circ-opt-spectral}
{\rm Opt}(A)=F_n {\rm diag}(F_n^*A F_n) F_n^*,
\end{equation}
where $A$ ia a generic complex square matrix and
the operator diag applied to any square matrix $X$ gives the diagonal matrix whose diagonal
entries coincide with those of $X$. Moreover, if $A=T_n(f)$ then
\begin{equation}\label{circ-opt}
{\rm Opt}(A)=C_n(f)={\rm circ}(a),\ \ \
\textstyle a_i=\frac{1}{n}\,((n-i)f_i+i f_{i-n})\qquad
(i=0,\dots,\lfloor n-1 \rfloor).
\end{equation}

In the nest proposition we discuss the spectral properties of these matrix approximations
by focusing on the relationships with the related approximation of the symbol.

\begin{proposition}\label{prop:general}
Let $f\in L^1(Q)$, $Q=(-\pi,\pi)$, and let us consider $N_n(f)$ and $C_n(f)$ be
the Strang and optimal approximations of $T_n(f)$, respectively. Then the following
facts hold:
\begin{enumerate}
\item The Strang preconditioner $N_n(f)$ has eigenvalues ${\cal F}_{n'}[f](x_j^{(n)})$,
$j=0,\ldots,n-1$, where $n'=[n/2]-1$, and ${\cal F}_{q}[f]$ is the Fourier sum of degree $q$ of $f$
(see {\rm \cite{CN}}).
\item In the general case where $f\in L^1(Q)$ and it is not smooth,
anything can happen: $N_n(f)$ definitely singular or indefinite
even if $T_n(f)$ is positive definite for every $n$,
$N_n(f)$ collectively unbounded even if $\|T_n(f)\|\le \|f\|_\infty$
for every $n$, $\{N_n(f)\}$ clustered at infinity even if $\{T_n(f)\}\sim_\sigma (f,Q)$.
\item If $f$ belongs to the Dini-Lipschitz class and is $2\pi$-periodic, then the eigenvalues of $N_n(f)$
will reconstruct $f$ in uniform norm.
\item The optimal preconditioner $C_n(f)={\rm Opt}(T_n(f))$ has eigenvalues ${\cal C}_{n-1}[f](x_j^{(n)})$,
$j=0,\ldots,n-1$, where ${\cal C}_{q}[f]={1\over q+1}\sum {\cal F}_{j=0}^{q}[f]$
is the Cesaro sum of degree $q$ of $f$ (see {\rm \cite{Skoro}}).
\item If $f$ is continuous and $2\pi$-periodic, then the eigenvalues of
$C_n(f)$ will reconstruct $f$ in uniform norm.
\item $\|{\rm Opt}(A)\|_*\le \|A\|_*$ for every unitarily invariant norm and in particular
$\|C_n(f)\|_p\le \|T_n(f)\|_p$ for every Schatten $p$ norm, $p\ge 1$ (see Theorem 2.1, item 6, in
{\rm \cite{fabio-lama}}).
\item If $f$ is $L^\infty(Q)$, then $\|C_n(f)\|\le \|f\|_\infty$ and, if $f\in L^p(Q)$
then $\|C_n(f)\|_p^p\le {n \over 2\pi} \int_Q |f(s)|^p ds$.
\item $\{C_n(f)\}$ distributes as $(f,Q)$ both in the sense of the eigenvalues and singular values.
\item With the notation of Proposition \ref{prop:sec:3-1},
$\{{\rm Opt}(\tilde M_n[\tilde \phi])\}$ distributes as $(\tilde f,Q)$ both in the sense of the
eigenvalues and singular values with
$\tilde f(s)=\pi\phi(\cos(s))w(\cos(s))\sin(s)$.
\end{enumerate}
\end{proposition}
{\bf Proof.}
{\bf Items 1., 4., 6.} can be found in the relevant literature.  {\bf Item 3.} is a direct consequence of
the fact that the Lebesgue constant of the Fourier sum is asymptotic (up to a multiplicative constant)
to $\log(n)$ and therefore the Fourier sum has to converge to $f$ since the modulus of continuity
of $f$ satisfies $\omega_f(1/n)=o(1/\log(n))$ for every $f$ in the Dini-Lipschitz class. {\bf Item 2.}
is a nice application of known facts. The example of Du Bois-Raymond is a nonnegative function
$f\in L^\infty(Q)$ with unbounded, highly oscillating Fourier sum (see e.g. \cite{bhatia-fourier}).
Clearly the matrix $N_n(f)$ is unbounded
and definitely indefinite while $T_n(f)$ is positive definite and uniformly bounded in spectral norm
by $\|f\|_\infty$ (for the Toeplitz part see e.g. \cite{jia} where also the tools for proving item 6
of Theorem 2.1 in {\rm \cite{fabio-lama}} can be found). For finding an example where
$\{N_n(f)\}$ clustered at infinity even if $\{T_n(f)\}\sim_\sigma (f,Q)$, it is enough to use the example
of Kolmogorov (see e.g. \cite{bhatia-fourier}): the function belongs to $L^1(Q)$, but it is not in $L^2(Q)$ and
has a Fourier sum diverging everywhere so that the eigenvalues of $N_n(f)$ collectively explode, but thanks
to Theorem \ref{teoszego-tyr} it is still true that $\{T_n(f)\}\sim_\sigma (f,Q)$.
{\bf Item 5.} is trivial since (thanks e.g. to the beautiful theory by Korovkin) it is well known
that the Cesaro sum of any continuous function converges uniformly to $f$. By \cite{jia}
we know that $\|T_n(f)\|\le \|f\|_\infty$ whenever $f\in L^\infty(Q)$ and
$\|T_n(f)\|_p^p\le {n \over 2\pi} \int_Q |f(s)|^p ds$ whenever $f\in L^p(Q)$ with $p\ge 1$:
as a consequence, {\bf Item 7.} follows from {\bf Item 6}.

Concerning {\bf Item 8.} we remark it has been
proved that for every $f\in L^1(Q)$, $\{T_n(f)\}\sim_\sigma (f,Q)$ and
$\{T_n(f)\}\sim_\lambda (f,Q)$ if $f$ is real-valued (see \cite{korotest}). We then need only to prove that the distribution results
stands for the eigenvalues as well even for complex-valued symbols (notice that the latter is not trivial
since it does not hold in general in the Toeplitz case as observed by Morrison in \cite{morrison}).

We want to prove that
\[
\lim_{n\rightarrow \infty}\Sigma_{\lambda}(F,C_n(f))=
{1\over 2\pi}\int_Q F(f(s))\, ds
\]
for every $f\in L^1(Q)$, for every $F$ continuous with bounded support in $\bf C$. First we observe
that the claim can be reduced to the case of $F$ Lipschitz continuous with bounded support in $\bf C$.
In fact for every $G$ continuous with bounded support in $\bf C$, for every $\epsilon>0$, we can find $G_\epsilon$
Lipschitz continuous with bounded support such that $|G(z)-G_\epsilon(z)|<\epsilon$ for every $z\in \bf C$
(notice that in general we cannot take $G_\epsilon$ polynomial due to the obstruction given by the
Mergelyan theorem (for a proof see \cite{Ru})).

Now by {\bf Item 5.} the claim is already proven if $f$ is continuous and $2\pi$-periodic (notice that in the
Toeplitz case this is again false in general with elementary polynomial examples).
Therefore for every $f\in L^1(Q)$, for every $\epsilon>0$, we consider
$f_\epsilon$ continuous and $2\pi$-periodic such that $\|f-f_\epsilon\|_{L^1(Q)}\le 2\pi\epsilon$ so that
\[
\left|{1\over 2\pi}\int_Q F(f(s))\, ds-{1\over 2\pi}\int_Q F(f_\epsilon(s))\, ds\right|
\le {1\over 2\pi}\int_Q |F(f(s))-F(f_\epsilon(s))|\, ds \le M\epsilon
\]
with $M$ being the Lipschitz constant of $F$.
Moreover, by the same argument we have,
\[
\left|\Sigma_{\lambda}(F,C_n(f))-\Sigma_{\lambda}(F,C_n(f_\epsilon))\right|
\le {1\over n}\sum_{j=1}^n |F(\lambda_j(C_n(f)))-F(\lambda_j(C_n(f_\epsilon)))|
\le M {1\over n}\sum_{j=1}^n |\lambda_j(C_n(f))-\lambda_j(C_n(f_\epsilon))|
\]
and, since the circulants form an algebra and the operator $C_n(\cdot)$ is linear, we have
\[
\left|\Sigma_{\lambda}(F,C_n(f))-\Sigma_{\lambda}(F,C_n(f_\epsilon))\right|
\le M {1\over n}\sum_{j=1}^n |\lambda_j(C_n(f-f_\epsilon))|.
\]
But the singular values of any circulant matrix are exactly the moduli of its eigenvalues since
every circulant is also normal. Therefore, by {\bf Item 7.}, we have
\[
\left|\Sigma_{\lambda}(F,C_n(f))-\Sigma_{\lambda}(F,C_n(f_\epsilon))\right|
\le {M\over n}\|C_n(f-f_\epsilon)\|_1\le M\epsilon
\]
and the proof is concluded since $\epsilon$ is arbitrary.

We conclude with the proof of {\bf Item 9.}
By {\bf Item 5.} and {\bf Item 6.} of Proposition \ref{prop:sec:3-1},
we have $\tilde M_n[\tilde \phi]=T_n(\tilde f)+H_n(\tilde f)$
and $\tilde f(s)=\pi\phi(\cos(s))w(\cos(s))\sin(s)$.
Therefore by linearity of the operator Opt$(\cdot)$ we deduce
${\rm Opt}(\tilde M_n[\tilde \phi])={\rm Opt}(T_n(\tilde f))+{\rm Opt}(H_n(\tilde f))
=C_n(f)+{\rm Opt}(H_n(\tilde f))$.
Now, by the previous item, $\{C_n(f)\}$ distributes as $\tilde f$ over $Q$ both in the sense of the
eigenvalues and singular values. Moreover, by \cite{FasTilli}, $\|H_n(\tilde f)\|_1=o(n)$
and therefore by {\bf Item 6.} $\|{\rm Opt}(H_n(\tilde f))\|_1=o(n)$. Furthermore, from Lemma \ref{distr-schatten}
we deduce $\{{\rm Opt}(H_n(\tilde f))\}\sim_\sigma (0,Q)$,
$\{{\rm Opt}(\tilde M_n[\tilde \phi])\}\sim_\sigma (\tilde f,Q)$ and
$\{{\rm Opt}(H_n(\tilde f))\}\sim_\lambda (0,Q)$,
$\{{\rm Opt}(\tilde M_n[\tilde \phi])\}\sim_\lambda (\tilde f,Q)$ if $\phi$ is real-valued.
Finally, for the complex-valued case when considering the distribution in the eigenvalue sense, the proof is as
in the preceding item.
\ \hfill $\bullet$
\ \\

\begin{oss}\label{carlesson}\rm
{\bf Item 1.} in the above proposition has an interesting consequence. Take $f\in L^\infty(Q)$
and consider $N_n(f)$. Since the entries of $N_n(f)$ contain exactly the same coefficients
as $T_{n'}(f)$ with every Fourier coefficient counted $2n'$ times it follows that
$\|N_n(f)\|_2^2={2n'\over 2\pi}\|{\cal F}_{n'}[f]\|_{L^2}^2\le {n\over 2\pi}\|f\|_{L^2}^2
\le n\|f\|_\infty^2$. Therefore, by the spectral decomposition of $N_n(f)$ in {\bf Item 1.},
it follows
\[
\|N_n(f)\|_2^2=\sum_{j=0}^{n-1}|{\cal F}_{n'}[f](x_j^{(n)})|^2\le n \|f\|_\infty^2.
\]
Consequently, the cardinality of the set of indices $j$ such that
${\cal F}_{n'}[f](x_j^{(n)})$ is unbounded as $n$ tends to infinity has to be $o(n)$ and the infinity
norm of ${\cal F}_{n'}[f]$ over the grid-sequence
$\{x_j^{(n)}\}_{n}$ is at most $O(\sqrt n)$.
This means that the set of grid points in which the Fourier sum can diverge is negligible and
more precisely its cardinality is $o(n)$. Taking into account the possible maximal growth of a polynomial
of degree $n'=[n/2]-1$, it follows that the set where the Fourier sum can diverge in $[-\pi,\pi]$
has to be of zero Lebesgue measure and this is a linear algebra version of a Carlesson-type result
(see e.g. \cite{bhatia-fourier}).
\end{oss}

\begin{oss}\label{circulants-vs-rest-of-the-world}\rm
In \cite{morrison}, the author observed that Toeplitz sequences are unable to reconstruct $f$, in general,
if $f$ is complex-valued. Tilli gave a precise answer by characterizing the cases where this reconstruction
is just impossible. In {\bf Item 8.}, we proved that a special circulant approximation of $T_n(f)$
is indeed able to reconstruct the symbol $f$ in the maximal generality that is for $f\in L^1(Q)$:
moreover, by {\bf Item 5.}, if $f$ is also continuous and $2\pi$-periodic, the reconstruction can be
performed in a strong sense i.e. in uniform norm. This
is confirmation of the great stability of the considered approximation which has two reasons: the first is the
normality of circulants (in contrast with Toeplitz matrices generated by complex-valued symbol
which can be of maximal non-normality as any Jordan block), the second is the stability of the Frobenius
optimal approximation which has to be related to the stability of Linear Positive Operators.
What we will discuss in the next subsection is interesting, because it shows that, under mild assumptions,
the problem of the identification and reconstruction of the multiplier $\phi$ can be reduced also to
the Frobenius optimal circulant approximation of a Toeplitz matrix generated by a $L^1(Q)$ symbol:
the theoretical basis relies on Proposition \ref{prop:general} and especially {\bf Items 5.},
{\bf 8.}, and {\bf 9}.
\end{oss}

\subsubsection{Circulant based algorithms}\label{subs:3-1-2}

We only suppose to know the coefficients of $M_n[\phi]$ and the weight $w$ with the related
matrix $L_n$ and unknown $\phi$ (the case where $\phi$ is known with unknown weight $w$ leads to
a different problem).
The algorithm is heavily related to the analysis in Proposition \ref{prop:sec:3-1} and in
Proposition \ref{prop:general}. It can be roughly sketched as follows:

\begin{enumerate}
\item Form $M_n[\phi]$ and from $L_n$ (known when the weight $w$ is known), compute
$X_n=\tilde M_n[\tilde \phi]$ with $\tilde \phi(x)=\phi(x) w(x) \sqrt{1-x^2}$;
\item compute $C_n$ the Frobenius optimal approximation of $X_n$;
\item compute the eigenvalues of $C_n$ by FFT (storing also the index of the related eigenvectors);
\item reconstruct the function $\tilde f(s)=\pi\phi(\cos(s))w(\cos(s))\sin(s)$ and therefore
dividing by $\pi w(\cos(s)) \sin(s)$ reconstruct $f(s)=\phi(\cos(s))$, $s\in Q$ i.e. $\phi(x)$, $x\in [-1,1]$.
\end{enumerate}

Indeed the correctness of the above procedure is based on the last item of Proposition \ref{prop:general},
since $X_n=\tilde M_n[\tilde \phi])$ and $C_n={\rm Opt}(\tilde M_n[\tilde \phi])$ (see also {\bf Item 5.}
and {\bf Item 8.}).

We observe that the matrix $X_n$ in view of (\ref{toep-hankel}) contains an Hankel part which
represents a disturbance. Therefore, also in order to exploit the computationally convenient
formula (\ref{circ-opt}), we can eliminate this part. The argument is a trivial application of
the Riemann-Lebesgue Lemma (see \cite{Ru}): indeed, instead of
$X_n=\tilde M_n[\tilde \phi]=T_n(\tilde f)+H_n(\tilde f)$ we would like
to consider the matrix $T_n(\tilde f)$ only. Unfortunately, the matrix $T_n(\tilde f)$ is unknown
(only the entries of the whole matrix $X_n$ are available) and we
will approximate it by the Toeplitz matrix $\tilde T_n$ constructed according to the following idea.
We have $(X_n)_{n,n}=(T_n(\tilde f))_{n,n}+(H_n(\tilde f))_{n,n}=f_0+f_{2n}\approx f_0$ since, by the
Riemann-Lebesgue Lemma, $f_{2n}$ is infinitesimal: we set $(\tilde T_n)_{j,j}=f_0+f_{2n}$, $j=1,\ldots,n$.
We observe that $(X_n)_{n-1,n}+(X_n)_{n,n-1}=(T_n(\tilde f))_{n-1,n}+(T_n(\tilde f))_{n,n-1}+
(H_n(\tilde f))_{n-1,n}+(H_n(\tilde f))_{n,n-1}=f_{-1}+f_1+2f_{2n-1}$. Since $\tilde f$ is even we
have $f_{-j}=f_j$ for all $j\in \bf Z$ and therefore $((X_n)_{n-1,n}+(X_n)_{n,n-1})/2=f_1+
f_{2n-1}\approx f_1$ since, by the
Riemann-Lebesgue Lemma, $f_{2n-1}$ is infinitesimal: we set $(\tilde T_n)_{j,j-1}=(\tilde T_n)_{j-1,j}
f_1+f_{2n-1}$, $j=2,\ldots,n$. We proceed by considering
$(X_n)_{n-2,n}+(X_n)_{n-1,n-1}+(X_n)_{n,n-2}=f_{-2}+f_0+f_2+3f_{2n-1}$. Now we already computed
an approximation of $f_0$ and we know that $f_{-2}=f_2$. Therefore we can compute $f_2$ within an infinitesimal
error since
\[
(\tilde T_n)_{j,j-2}=(\tilde T_n)_{j-2,j}
=\left[((X_n)_{n-2,n}+(X_n)_{n-1,n-1}+(X_n)_{n,n-2})-(X_n)_{n,n}\right]/2=
f_2+(3f_{2n-1}-f_{2n})/2\approx f_2,
\]
$j=3,\ldots,n$. The procedure can be continued in a similar way by obtaining every entry of $T_n(\tilde f)$
i.e. every Fourier coefficient $f_j$, $|j|\le n-1$, within an infinitesimal approximation error.

The new algorithm can be written as follows where the  third step is obtained by the previous reasoning.

\begin{enumerate}
\item Form $M_n[\phi]$ and from $L_n$ (known when the weight $w$ is known), compute
$X_n=\tilde M_n[\tilde \phi]$ with $\tilde \phi(x)=\phi(x) w(x) \sqrt{1-x^2}$;
\item compute $\tilde T_n$ approximation of $T_n(\tilde f)$;
\item compute $C_n$ the Frobenius optimal approximation of $\tilde T_n$;
\item compute the eigenvalues of $C_n$ by FFT (storing also the index of the related eigenvectors);
\item reconstruct the function $\tilde f(s)=\pi\phi(\cos(s))w(\cos(s))\sin(s)$ and therefore
dividing by $\pi w(\cos(s)) \sin(s)$ reconstruct $f(s)=\phi(\cos(s))$, $s\in Q$ i.e. $\phi(x)$, $x\in [-1,1]$.
\end{enumerate}

The correctness of the algorithm depends entirely on the fact that $\{\tilde T_n\}$ distributes as $(\tilde f,Q)$
in the sense of the singular values. We know $\{T_n(\tilde f)\}\sim_\sigma (\tilde f,Q)$ and
$|(\tilde T_n-T_n(\tilde f))_{j,k}|$ tends to zero for every $(j,k)$
(more precisely we have $|(\tilde T_n-T_n(\tilde f))_{j,k}|=O(f_n)$). Unfortunately, the second relation
does not imply that $\|\tilde T_n-T_n(\tilde f)\|_p=o(n^{1/p})$ for some $p\in [1,\infty]$
and therefore we are not allowed to use Lemma \ref{distr-schatten}.
In fact, let us consider the following example. Assume that
\[
\tilde T_n-T_n(\tilde f)=\epsilon_n \sqrt n F_n,\ \ \ \epsilon_n>0\ \ \rm infinitesimal.
\]
Then every entry has modulus $\epsilon_n$ but every eigenvalue has modulus equal to $\sqrt n\epsilon_n$
and therefore $\{\tilde T_n-T_n(\tilde f)\}$ distributes as $(0,Q)$ only if $\epsilon_n=o(n^{-1/2})$.
Indeed, defining $\epsilon_n$ the maximal values of $|(\tilde T_n-T_n(\tilde f))_{j,k}|$,
we have
\[
\|\tilde T_n-T_n(\tilde f)\|_2^2=\sum_{j,k=0}^{n-1}|(\tilde T_n-T_n(\tilde f))_{j,k}|^2 \le
\sum_{j,k=0}^{n-1}\epsilon_n^2=\epsilon_n^2 n^2.
\]
Therefore, for $p=2$ we obtain $\|\tilde T_n-T_n(\tilde f)\|_p=o(n^{1/2})$ if $f_n=o(n^{-1/2})$.
As a consequence, in order to use the second
algorithm, we should have more information on the symbol $\tilde f$: for instance if $\tilde f$ is
$2\pi$-periodic and $k$-times continuously differentiable, $k\ge 1$, then $f_n=o(n^{-k})$ and therefore
$\epsilon_n=o(n^{-k})$ so that we can use safely the second algorithm and, in addition, if $k$ is large then
$\tilde T_n$ is a very good approximation of $T_n(\tilde f)$. In conclusion, for a smooth symbol
i.e. for large $k$, the reconstruction provided by the latter algorithm could be better than the one
given by the first. Of course, one should have this information a priori or, possibly, one should
use the first algorithm for obtaining a guess: if the result looks like a smooth function (under the assumptions
of {\bf Item 5.} in Proposition \ref{prop:general}, we recall that the approximation is in uniform norm), then
the second algorithm could be employed for improving the quality of the reconstruction.

\subsection{Generalization: the multivariate separable case}\label{subs:3-2}

Let $\phi$ be a bounded function defined on $[-1,1]^d$ and let us consider the multiplication
operator $M[\phi]:L_w^2([-1,1]^d)\rightarrow L_w^2([-1,1]^d)$ defined as $M[\phi](f)=\phi f$
with $w(x)=w_1(x_1)w_2(x_2)\cdots w_d(x_d)$, $d$-variate weight with $w_j$ standard univariate weight function.
As in the single-variate case the spectrum coincides with the essential range of $\phi$.
Consider to have a finite number of coefficients
\begin{equation}\label{form-m-phi-d}
M_n[\phi]=\left(\langle M[\phi]e_j,e_i\rangle\right)_{i,j=0}^{n-{\bf 1}^T},
\end{equation}
$n=(n_1,\ldots,n_d)$, $i=(i_1,\ldots,i_d)$, $j=(j_1,\ldots,j_d)$, ${\bf 1}$ vector of all ones
as in Subsection \ref{subs:2-1}, with $\{e_j\}$ denoting an orthonormal basis of $L_w^2([-1,1]^d)$ defined by
$e_j(x)=e_{j_1}(x_1)\cdots e_{j_d}(x_d)$ with  $\{e_{j_k}\}$ orthonormal basis of $L_{w_k}^2([-1,1])$,
$k=1,\ldots,d$. The question is again the reconstruction of the multiplier $\phi$ from the spectra
of $M_n[\phi]$ in the sense discussed in Subsection \ref{sec:2-1}.

Take the case of the $d$-level Chebyshev weight of first kind $w(x)=\prod_{k=1}^d (1-x_{j_k}^2)^{-1/2}$
and of its basis $e_j(x)=e_{j_1}(x_1)\cdots e_{j_d}(x_d)$, $e_{j_k}(x_k)=\cos(j_k\arccos(x_k))$,
$k=1,\ldots,d$. Then with the usual change of variable $x_k=\cos(s_k)$, $k=1,\ldots,d$, we have
\[
\left(M_n[\phi]\right)_{i,j}=\langle M[\phi]e_j,e_i\rangle=
\int_{Q^d}\phi(\cos(s_1),\ldots,\cos(s_d))\prod_{k=1}^d\cos(j_ks_k)\cos(i_ks_k)\, ds,\ \ Q=(-\pi,\pi).
\]
As a consequence, the matrix $M_n[\phi]$ can be expressed in terms of the $d$-indexed
Fourier coefficients $f_j$ of the function $f(s)=(\pi)^d \phi(\cos(s_1),\ldots,\cos(s_d))$ and then,
in $d$-index notation and taking into account that $f(s)$ is even with respect to every variable
$s_j$, we find
\[
\left(M_n[\phi]\right)_{i,j}=f_{|i-j|}+f_{|i+j|}
\]
and therefore
\begin{equation}\label{toep-hankel-d}
M_n[\phi]=T_n(f)+H_n(f).
\end{equation}
Once we arrive here, the rest is a straightforward generalization of the univariate case since
the result on Hankel matrices (see \cite{FasTilli}) are directly stated in an arbitrary number of
dimensions. Theorem \ref{teoszego-tyr} is in $d$ dimensions and the same applies to the results on the GLT class
whose definition is inherently $d$-dimensional (see Subsection \ref{subs:2-2} and \cite{multiloc}).
Furthermore, formulae (\ref{circ-opt-spectral})--(\ref{circ-opt}) stand unchanged
($d$-indices in place on simple indices, $F_n=F_{n_1}\otimes \cdots \otimes F_{n_d}$,
$n$ in the denominator of (\ref{circ-opt}) replaced by $N(n)$)
and Proposition \ref{prop:general} is again unchanged. Therefore also the algorithms can be described
verbatim and therefore the optimal circulant approximation of $M_n[\phi]\approx T_n(f)$,
$\tilde T_n\approx T_n(f)$ will reconstruct with infinitesimal error the Cesaro sum of the function
$f$ and therefore of $\phi$.


\subsection{Generalization: the multivariate separable matrix-valued case}\label{subs:3-3}

Let $\phi$ be a bounded function defined on $[-1,1]^d$ and having values in the space
${\bf C}^{p\times q}$ and let us consider the multiplication
operator $M[\phi]:L_w^2([-1,1]^d,{\bf C}^{q\times r})\rightarrow
L_w^2([-1,1]^d,{\bf C}^{p\times r})$ defined as $M[\phi](f)=\phi f$ being
${\bf C}^{p\times r}$ with $f$ being ${\bf C}^{q\times r}$ and
with $w(x)=w_1(x_1)w_2(x_2)\cdots w_d(x_d)$ as in the previous subsection.
In the present matrix-valued setting, it is less obvious to refer to the spectrum of the
continuous operator and this is true in the discrete as well since the resulting sections
are not square matrices. However we can give again a meaning passing to the ``absolute value''
of the operator (see \cite{Bha})
i.e. the square root of the adjoint times the operator itself. In the discrete we are talking
of the singular values and in the operator case we are talking of the singular values of the
multiplier $\phi$. In any case, as we will see in the rest of the derivation, we will able
to reconstruct $\phi$ (or better its Cesaro sum) and therefore its singular values.

It should be observed that the present multiplication operator can be written as a vector whose entries are
sum of scalar multiplication operators. More precisely we have
\[
\phi f=\left(\sum_{t=1}^{q}\phi_{s,t}f_{t,z}\right)_{s=1,z=1}^{p,r}
\]
with $\phi_{s,t}f_{t,z}$ defining a scalar multiplication operator on $L_w^2([-1,1]^d)$.
Therefore we can represent $M[\phi]$ as
\[
\sum_{s=1}^{p}\sum_{t=1}^{q}M[\phi_{s,t}E(s,t)]
\]
where $E(s,t)$ denotes the $p\times q$ matrix being 1 at position $(s,t)$ and zero
otherwise: notice that $\{E(s,t)\}$ forms the canonical basis ${\bf C}^{p\times q}$.
Therefore, if we consider a finite number of coefficients we have
\begin{eqnarray*}
M_n[\phi] & = & \left(\langle M[\phi]e_j,e_i\rangle
\right)_{i,j=0}^{n-{\bf 1}^T} \\
 & = & \left(\langle \sum_{s=1}^{p}\sum_{t=1}^{q}M[\phi_{s,t}E(s,t)]e_j,e_i\rangle
\right)_{i,j=0}^{n-{\bf 1}^T} \\
& = & \sum_{s=1}^{p}\sum_{t=1}^{q}\left(\langle M[\phi_{s,t}E(s,t)]e_j,e_i\rangle
\right)_{i,j=0}^{n-{\bf 1}^T},
\end{eqnarray*}
$n=(n_1,\ldots,n_d)$, $i=(i_1,\ldots,i_d)$, $j=(j_1,\ldots,j_d)$, ${\bf 1}$ vector of all ones,
with $\{e_j\}$ denoting an orthonormal basis of $L_w^2([-1,1]^d)$ as in the latter
section.
In other words $M_n[\phi]$ is a multilevel block matrix where the size of each block is dictated
by the multiplier $\phi$ and more specifically we can write
\begin{equation}\label{form-m-phi-block}
\left(M_n[\phi]\right)_{s,t}=\left[\left(\langle M[\phi_{s,t}]e_j,e_i)\rangle\right)_{s,t}
\right]_{s=1,t=1}^{p,q}=
\left[ \begin{array}{ccc}
    \left(M_n[\phi_{1,1}]\right)_{i,j} &  \cdots & \left(M_n[\phi_{1,q}]\right)_{i,j} \\
    \vdots & \ddots  & \vdots  \\
   \left(M_n[\phi_{p,1}]\right)_{i,j} &  \cdots & \left(M_n[\phi_{p,q}]\right)_{i,j}
     \end{array} \right].
\end{equation}

The above block expression makes clear that the reconstruction of every single entry
$\phi_{s,t}$ can be done exactly as in the scalar-valued case. More precisely the next scheme can be
followed.
\begin{itemize}
\item The reconstruction of the entry $\phi_{s,t}$ can be done via the same algorithms proposed in Subsection
\ref{subs:3-1-2}, by extracting from $M_n[\phi]$ only the entries of $M_n[\phi_{s,t}]$ according to
(\ref{form-m-phi-block}).
\item If one is interested in a global reconstruction e.g. of the singular values of the multiplier $\phi$,
then some appropriate new tools have to be introduced: this issue is briefly considered in the next section.
\end{itemize}

\subsubsection{Toeplitz sequences generated by matrix-valued sumbols}\label{subs:3-3-1}

Let $f$ be a $d$ variate $p\times q$ matrix-valued (Lebesgue) integrable function,
defined over the hypercube $Q^d$, with $Q=(-\pi,\pi)$ and $d\ge 1$. Here the Lebesgue integrability
means that every entry of the symbol is a standard complex-valued $L^1$ function.
With respect to the notion of matrix-valued symbol, we observe that the definition of the coefficients
in (\ref{defcoeff}) is formally identical: the only obvious difference is that every $f_j$ will be a matrix
of size $p\times q$. However, the formulae in (\ref{distribution:sv-eig}) and consequently Theorem
\ref{teoszego-tyr} do not make sense since the eigenvalues and
singular values are still scalar while the function to be integrated in the right hand-side is
$p\times q$ matrix-valued. On the other hand, a generalization of that results exists and is quite natural.
We write that $\{A_n\}\sim_{\sigma} (f,K,\mu)$
and $\{A_n\}\sim_{\lambda} (f,K,\mu)$, if
\begin{eqnarray}\label{distribution:sv-eig-block}
\lim_{n\rightarrow \infty}\Sigma_{\sigma}(F,A_n)& = & {1\over \mu(K)}\int_K {1\over l}{\rm tr}[F(|f|)]d\mu, \\
\nonumber
\lim_{n\rightarrow \infty}\Sigma_{\lambda}(F,A_n) & = &
{1\over \mu(K)}\int_K {1\over l}{\rm tr}[F(f)]d\mu,\ \ l=p=q,\ f\  \hbox{\rm Hermitian-valued},
\end{eqnarray}
respectively, with $l$ being the minimum between $p$ and $q$, $|f|=(f^H f)^{1/2}$ and ${\rm tr}[g]=\sum_j \lambda_j(g)$, $\lambda_j(g)$,
$j=1,\ldots,l$, being the eigenvalues of $g$. With these notations, we have that any Toeplitz sequence
$\{T_n(f)\}$ with matrix-valued $f\in L^1(Q^d)$ (which is equivalent to require that maximal singular
value of $f$ is in $L^1(Q^d)$) is such that (\ref{distribution:sv-eig-block}) holds with
$\mu$ being the Lebesgue measure, $K=Q^d$, and $f$ being the symbol
(see \cite{tilliLA,Tillinota,Smarko1}).
Notice that if $g=|f|$ then ${\rm tr}[g]=\sum_j \sigma_j(f)$ and therefore (\ref{distribution:sv-eig-block})
represents a natural generalization of \rif{distribution:sv-eig}.

The question is about the reconstruction of the multiplier $\phi$ from the finite section
$M_n[\phi]$.

Take the case of the $d$-level Chebyshev weight of first kind $w(x)=\prod_{k=1}^d (1-x_{j_k}^2)^{-1/2}$
and of its basis $e_j(x)=e_{j_1}(x_1)\cdots e_{j_d}(x_d)$, $e_{j_k}(x_k)=\cos(j_k\arccos(x_k))$,
$k=1,\ldots,d$. Then with the usual change of variable $x_k=\cos(s_k)$, $k=1,\ldots,d$, we have
\[
\left(M_n[\phi]\right)_{i,j}=\langle M[\phi]e_j,e_i\rangle=
\int_{Q^d}\phi(\cos(s_1),\ldots,\cos(s_d))\prod_{k=1}^d\cos(j_ks_k)\cos(i_ks_k)\, ds,\ \ Q=(-\pi,\pi).
\]
As a consequence, the matrix $M_n[\phi]$ can be expressed in terms of the $d$-indexed
Fourier coefficients $f_j$ of the function $f(s)=(\pi)^d \phi(\cos(s_1),\ldots,\cos(s_d))$ and then,
in $d$-index notation, we find
\[
\left(M_n[\phi]\right)_{i,j}=f_{|i-j|}+f_{|i+j|}.
\]
Taking into account that $f(s)$ is even we directly see that
\begin{equation}\label{toep-hankel-block}
M_n[\phi]=T_n(f)+H_n(f).
\end{equation}
Once we arrive here, the rest is a generalization of the
multivariate case since the result on Hankel matrices (see
\cite{FasTilli}) are directly stated in an arbitrary number of
dimensions with blocks of fixed dimension. Theorem
\ref{teoszego-tyr} is replaced by the block relation
(\ref{distribution:sv-eig-block}) and GLT class has a natural
block generalization. Furthermore, the formula (\ref{circ-opt})
stands unchanged ($d$-indices and block coefficients in place on
simple indices and scalar coefficients) and Proposition
\ref{prop:general} is again unchanged. Therefore also the
algorithms can be described verbatim and therefore the optimal
circulant approximation of $M_n[\phi]\approx T_n(f)$, $\tilde
T_n\approx T_n(f)$ will reconstruct with infinitesimal error the
Cesaro sum of the function $f$ and therefore of $\phi$.
Furthermore, by using (\ref{distribution:sv-eig-block}), it is
possible to reconstruct information on the singular values of
$f$ and then of $\phi$. More precisely, given $s\in \bf R$,
$s\ge 0$, there is constructive test, analogous to those in {\bf
e)} and {\bf f)} of Theorem \ref{th-relations-def}, that, starting
from the singular values of the matrix $M_n[\phi]\approx T_n(f)$
(or $\tilde T_n\approx T_n(f)$), tells one if $s$ belongs to the
union of the (essential) ranges of the singular values of $f$ (and
therefore of $\phi$): for details and numerical experiments in
this block setting see \cite{korotest}.


\section{Conclusions}\label{sec:4}

In this paper we have considered the reconstruction of the scalar-valued/multivariate/block-valued multipliers
of proper multiplication operators through structured linear algebra tools.
Further generalizations could be considered as a general compact domain or non-separable weight functions.
However, we think that the results presented in this note clearly show the utility of purely asymptotic
linear algebra tools in the considered type of problems in approximation theory.



\begin{thebibliography}{99}

\bibitem{Bha}
 R. Bhatia,
 {\em Matrix Analysis}.
  Springer-Verlag, New York, 1997.

\bibitem{bhatia-fourier}
 R. Bhatia,
 {\em Fourier Series}.
 AMS, Providence, 2004.

\bibitem{BG}
A. B\"ottcher and S. Grudsky,
``On the condition numbers of large semi-definite Toeplitz
matrices'',
{\it Linear Algebra Appl.}, {\bf 279} (1998), pp. 285--301.

\bibitem{BS}
A. B\"ottcher and B. Silbermann,
{\em Introduction to Large Truncated Toeplitz Matrices}.
Springer-Verlag, New York, 1999.

\bibitem{CN}
 R.H. Chan and M. Ng,
  ``Conjugate gradient methods for Toeplitz systems'',
 {\it SIAM Rev.}, {\bf 38} (1996), pp. 427--482.

\bibitem{Da}
P. Davis,
{\em Circulant Matrices.}
John Wiley and Sons, New York, 1979.

\bibitem{fabio-lama}
F. Di Benedetto and S. Serra-Capizzano,
``Optimal multilevel matrix algebra operators'',
{\it Linear Multilin. Algebra}, {\bf 48} (2000), pp. 35--66.


\bibitem{FasTilli}
D. Fasino and P. Tilli,
``Spectral clustering properties of block multilevel Hankel matrices'',
{\it Linear Algebra Appl.}, {\bf 306} (2000), pp. 155--163.

\bibitem{Salmeria}
L. Golinskii and S. Serra-Capizzano, ``The asymptotic spectrum of
non symmetrically perturbed symmetric Jacobi matrix sequences'',
{
http://arxiv.org/abs/math.SP/0512222}, (2005).

\bibitem{szego}
U.~Grenander and G.~{Szeg\H o}.
{\em Toeplitz Forms and their Applications}.
Second Edition, Chelsea, New York, 1984.

\bibitem{KS}
T. Kailath and A.H. Sayed,
``Displacement structures: theory and applications'',
{\it SIAM Rev.}, {\bf 37}  (1995), pp. 297--386.

\bibitem{ku-ser}
A.B.J. Kuijlaars and S. Serra-Capizzano,
``Asymptotic zero distribution of orthogonal polynomials with discontinuously
varying recurrence coefficients'', {\it J. Appr. Theory},
{\bf  113} (2001), pp. 142--155.

\bibitem{morrison}
K. Morrison,
``Spectral approximation of multiplication operators'',
{\it New York J. Math.},
{\bf 1} (1995), pp. 75--96.

\bibitem{Ru}
W. Rudin,
{\em Real and Complex Analysis}.
McGraw-Hill, New York, 1974.

\bibitem{Se1}
S. Serra-Capizzano,
``On the extreme eigenvalues of Hermitian (block) Toeplitz matrices'',
{\it Linear Algebra Appl.}, {\bf 270} (1998), pp. 109--129.



\bibitem{Smarko1}
S. Serra-Capizzano,
``Spectral and computational analysis of block Toeplitz matrices with
nonnegative definite generating functions'',
{\it BIT}, {\bf 39} (1999), pp. 152--175.

\bibitem{Skoro}
S. Serra-Capizzano,
``A Korovkin-type Theory for finite Toeplitz operators
via matrix algebras'',
{\it Numer. Math.}, {\bf 82-1} (1999), pp. 117--142.

\bibitem{korotest}
S. Serra-Capizzano
``Korovkin tests, approximation, and ergodic theory'',
{\it Math. Comput.},  {\bf 69} (2000), pp. 1533--1558.

\bibitem{algebra}
S. Serra-Capizzano,
``Distribution results on the algebra generated by Toeplitz
sequences: a finite dimensional approach'',
{\it Linear Algebra Appl.}, {\bf 328} (2001), pp. 121--130.

\bibitem{Spal}
S. Serra-Capizzano,
``Test functions, growth conditions and Toeplitz matrices'',
{\it Rend. Circolo Mat. Palermo}, Serie II, {\bf 68} (2002), pp. 791--795.

\bibitem{multiloc}
S. Serra-Capizzano,
``Generalized Locally Toeplitz sequences: spectral
analysis and applications to discretized Partial Differential equations'',
{\it Linear Algebra Appl.}, {\bf 366-1} (2003), pp. 371--402.

\bibitem{cluster}
S. Serra-Capizzano, D. Bertaccini, and G. Golub
{``How to deduce a proper eigenvalue cluster
from a proper singular value cluster in the non-normal case''},
{\it SIAM J. Matrix Anal. Appl.}, {\bf 27-1} (2005), pp. 82--86.

\bibitem{nota}
S. Serra-Capizzano,
``GLT sequences and orthogonal polynomials'',
{\it manuscript}, (2005).

\bibitem{sccm}
S. Serra-Capizzano,
{``GLT sequences as a Generalized Fourier Analysis and applications''},
{\it SCCM-5-7}, (2005), Stanford University.

\bibitem{jia}
S. Serra-Capizzano and P.~Tilli,
``On unitarily invariant norms of matrix
valued linear positive operators'',
{\it J. Inequalities Appl.}, {\bf 7-3} (2002), pp. 309--330.

\bibitem{szego-pol}
G.~{Szeg\H o}.
{\em Orthogonal Polynomials}.
AMS, fourth edition, 1978.

\bibitem{Tilliloc}
P. Tilli,
``Locally Toeplitz matrices: spectral theory and applications'',
{\it Linear Algebra Appl.}, {\bf 278} (1998), pp. 91--120.

\bibitem{Tillinota}
P. Tilli,
``A note on the spectral distribution of Toeplitz matrices'',
{\it Linear Multilin. Algebra}, {\bf 45} (1998), pp.~147--159.

\bibitem{tilliLA}
P.~Tilli,
``Singular values and eigenvalues of non-{Hermitian} block {Toeplitz} matrices''
{\it Linear Algebra Appl.}, {\bf 272} (1998), pp. 59--89.

\bibitem{tillicomplex}
   P. Tilli,
 ``Some results on complex Toeplitz eigenvalues'',
{\it J. Math. Anal. Appl.}, {\bf 239-2} (1999), pp. 390--401.

\bibitem{tyrtyshnikov}
 E. Tyrtyshnikov,
  ``A unifying approach to some old and new theorems on
    distribution and clustering'',
{\it Linear Algebra Appl.}, {\bf 232} (1996), pp. 1--43.

\bibitem{tyrtL1}
E. Tyrtyshnikov and N. Zamarashkin,
``Spectra of multilevel Toeplitz matrices: advanced theory via simple
matrix relationships'',
{\it Linear Algebra Appl.}, {\bf 270} (1998), pp. 15--27.

\bibitem{fft}
C.F. Van Loan, {\em Computational Frameworks for the Fast Fourier
Transform}. SIAM, Philadelphia, 1992.

\bibitem{Wil}
J. Wilkinson, {\em The Algebraic Eigenvalue Problem.}
Claredon Press, Oxford, 1965.

\bibitem{Zy}
A. Zygmund,
{\em Trigonometric Series.}
Cambridge University Press, Cambridge, 1959.

\end{thebibliography}
\end{document}